\documentclass{amsart}
\usepackage{amsthm}
\usepackage{amsmath}

\theoremstyle{plain}

\theoremstyle{plain}
\newtheorem{thm}{Theorem}[section]

\newtheorem{lem}[thm]{Lemma}

\newtheorem*{rem*}{Remark}

\numberwithin{equation}{section}

\email{inoue@math.sci.hokudai.ac.jp}
\email{nakano\_y@math.sci.hokudai.ac.jp}
\email{v.anh@fsc.qut.edu.au}
\thanks{This work is partially supported by the Australian 
 Research Council grant DP0345577.}

\begin{document}
\title[\null]{Binary market models with memory}
\author[\null]{Akihiko Inoue, Yumiharu Nakano and Vo Anh}
\address{Department of Mathematics \\
Faculty of Science \\
Hokkaido University \\
Sapporo 060-0810 \\
Japan}
\address{Department of Mathematics \\
Faculty of Science \\
Hokkaido University \\
Sapporo 060-0810 \\
Japan}
\address{School of Mathematical Sciences \\
Queensland University of Technology \\
GPO Box 2434, Brisbane, Queensland 4001 \\
Australia}
\begin{abstract}
We construct a binary market model with memory that approximates a
continuous-time market model driven by a Gaussian process equivalent to
Brownian motion. We give a sufficient conditions for the binary market to 
be
arbitrage-free. In a case when arbitrage opportunities exist, we present 
the
rate at which the arbitrage probability tends to zero as the number of
periods goes to infinity.
\end{abstract}

\maketitle


\section{Introduction}

\label{sec:1}

Let $T\in (0,\infty )$. We consider the stock price process $(S_{t})_{0\leq
t\leq T}$ that is governed by the stochastic differential equation 
\begin{equation}
dS_{t}=S_{t}(bdt+\sigma dY_{t})\qquad (0\leq t\leq T),  \label{eq:1.1}
\end{equation}
where $\sigma $ and the initial value $S_{0}$ are positive constants, and 
$%
b\in \mathbf{R}$. In the classical Black-Scholes model, Brownian motion is
used as the driving noise process $Y$, and the resulting price process $S$
becomes Markovian. In \cite{AI, AIK}, the following Gaussian process $%
(Y_t)_{0\leq t\leq T}$ with stationary increments is used instead as the
driving noise process $Y$ in (\ref{eq:1.1}): 
\begin{equation}
Y_{t}=B_{t}-\int_{0}^{t}\left\{ \int_{-\infty
}^{s}pe^{-(q+p)(s-u)}dB_{u}\right\} ds\qquad (0\leq t\leq T),  
\label{eq:1.2}
\end{equation}
where $p$ and $q$ are real constants such that 
\begin{equation*}
0<q<\infty ,\quad -q<p<\infty ,
\end{equation*}
and $(B_{t})_{t\in \mathbf{R}}$ is a one-dimensional Brownian motion 
defined
on a probability space $(\Omega ,\mathcal{F},P)$ satisfying $B_{0}=0$. The
parameters $p$ and $q$ describe the memory of $Y$, and the resulting stock
price process $S$ becomes non-Markovian. An empirical study on S\&P 500 
data
in \cite{AIP} shows that the model captures very well the memory effect 
when
the market is stable.

It should be noticed that (\ref{eq:1.2}) is not a semimartingale
representation of $Y$ with respect to the $P$-augmentation $(\mathcal{F}%
_{t})_{0\leq t\leq T}$ of the filtration generated by $(Y_{t})_{0\leq t\leq
T}$ since $(B_{t})$ is not $(\mathcal{F}_{t})$-adapted. However, by
innovation theory as described in Liptser and Shiryayev \cite{LS}, we can
show that $Y$ is actually an $(\mathcal{F}_{t})$-semimartingale 
(\cite[Theorem 3.1]{AI}). In fact, using the prediction theory for $Y$ which
is developed in \cite{AIK}, we see (\cite[Theorem 2.1]{INA}) that there
exists a one-dimensional Brownian motion $(W_{t})_{0\leq t\leq T}$, called
the \textit{innovation process}, satisfying 
\begin{equation*}
\sigma (W_{s}:0\leq s\leq t)=\sigma (Y_{s}:0\leq s\leq t)\qquad (0\leq 
t\leq
T),
\end{equation*}
and 
\begin{equation}
Y_{t}=W_{t}-\int_{0}^{t}\left\{ \int_{0}^{s}l(s,u)dW_{u}\right\} ds\qquad
(t\in \lbrack 0,T]),  \label{eq:1.3}
\end{equation}
where $l(t,s)$ is a Volterra kernel given explicitly by 
\begin{equation}
l(t,s)=pe^{-(p+q)(t-s)}\left\{ 
1-\frac{2pq}{(2q+p)^{2}e^{2qs}-p^{2}}\right\}
\quad (0\leq s\leq t\leq T).  \label{eq:1.4}
\end{equation}
Thus the process $Y$ has the virtue that it possesses the property of a
stationary increments process with memory and the simple semimartingale
representation (\ref{eq:1.3}) with (\ref{eq:1.4}) simultaneously. We know 
of
no other process with this kind of properties. The two properties of $Y$
become a great advantage, for example, in its parameter estimation (see 
\cite[Section 5]{INA}).

Several authors use fractional Brownian motion as the driving noise process 
(see, e.g., Comte and Renault \cite{CR}, Rogers 
\cite{R}, and Willinger et al.~\cite{WTT}). However this approach is not 
entirely 
satisfactory since fractional Brownian motion is not a semimartingale (Lin 
\cite{L} and Rogers \cite{R}), whence there exists no equivalent martingale
measure in the corresponding market. On the other hand, the market defined
by (\ref{eq:1.1}) with (\ref{eq:1.2}) or (\ref{eq:1.3}) and (\ref{eq:1.4})
is arbitrage-free and complete since the process $Y$ becomes a Brownian
motion under a suitable probability measure (see \cite[Section 3]{AI}).
Moreover, for this model, we can obtain explicit results such as the
solution to the expected logarithmic utility maximization from terminal
wealth (see \cite{AIK}).

As is well known, binary approximation of the Black-Scholes model plays a
very important role for the model in many ways. Sottinen \cite{S}
constructed a binary market model that approximates the market driven by
fractional Brownian motion, and investigated the arbitrage opportunities in
the binary model.

In this paper, we construct a binary market model with memory that
approximates the continuous-time market model driven by $Y$ in 
(\ref{eq:1.3}%
). However, rather than considering the special kernel $l(t,s)$ in (\ref
{eq:1.4}), we take a general bounded measurable Volterra kernel $l(t,s)$.
Since $l(t,s)$ given by (\ref{eq:1.4}) is bounded, the results thus 
obtained
apply to the special case (\ref{eq:1.4}). We remark that any centered
Gaussian process $Y=(Y_{t})_{0\leq t\leq T}$ that is equivalent to a
Brownian motion has a canonical representation of the form (\ref{eq:1.3})
with $l(t,s)$ satisfying square integrability (see Hida and Hitsuda 
\cite[Chapter VI]{HH}). Thus, in this paper, we consider a subclass
consisting of $Y$ for which $l(t,s)$ is bounded. As in \cite{S}, the key
feature to the construction of the approximating binary market is to prove 
a
Donsker-type theorem for the process $Y$ (Theorem \ref{thm:2.1}).

Unlike the market driven by fractional Brownian motion, the market driven 
by 
$Y$ in (\ref{eq:1.3}) is arbitrage-free (see, e.g., the proof of 
\cite[Theorem 3.3]{AI}). However, the approximating binary market may admit
arbitrage opportunities. We consider conditions for their existence or
non-existence. We also study the rate at which the aribtrage probability
tends to zero as the number of periods goes to infinity.

This paper is organized as follows. In Section 2, we prove a Donsker-type
theorem for the driving process $Y$ in (\ref{eq:1.3}) with bounded kernel 
$%
l(t,s)$. In Section 3, we consider a discrete-time approximation of the
stock price process $S$ in (\ref{eq:1.1}). As a special case, we obtain the
desired approximating binary model. In Section 4, we study arbitrage
opportunities in the binary model.


\section{A Donsker-type theorem }
\label{sec:2}

Let $T\in (0,\infty )$. In what follows, we write $C=C_{T}$ for positive
constants, depending on $T$, which may not be necessarily equal to each
other. Let $n$ be a positive integer. In Sections 2 and 3, we 
write
\[
\sum_{s\le t}X_s=\sum\nolimits_{i=1}^{\lfloor nt\rfloor}X_{\frac{i}{n}},\qquad 
\prod_{s\le t}X_s=\prod\nolimits_{i=1}^{\lfloor nt\rfloor}X_{\frac{i}{n}}.
\]

Let $l(t,s)$ be a bounded measurable function on $[0,T]\times [0,T]$ 
that vanishes 
whenever $s>t$. Let $W=(W_t)_{0\le t\le T}$ be a one-dimensional Brownian
motion on a probability space $(\Omega,\mathcal{F},P)$. We define the
process $Y=(Y_t)_{0\le t\le T}$ by (\ref{eq:1.3}).

We put, for $t, u\in [0,T]$, 
\begin{equation*}
z(t,u):= \int_{u}^{t}l(s,u)ds,\qquad y(t,u):=1-z(t,u). 
\end{equation*}
Then both $z(t,u)$ and $y(t,u)$ are bounded and continuous on $[0,T]\times
[0,T]$, and it holds that $Y_t=\int_{0}^{t}y(t,u)dW_u$ for $0\le t\le T$.
Let $C$ be a positive constant satisfying, for $(t_1,u), (t_2,u)\in
[0,T]\times [0,T]$, 
\begin{equation}
\vert z(t_1,u)-z(t_2,u)\vert =\vert y(t_1,u)-y(t_2,u)\vert \le C\vert
t_1-t_2\vert.  \label{eq:2.1}
\end{equation}

Let $\{\xi_i\}_{i=1}^{\infty}$ be a sequence of i.i.d.~random variables 
with 
$E[\xi_1]=0$ and $E[(\xi_1)^2]=1$. We also assume that 
\begin{equation}
E[(\xi_1)^4]<\infty.  \label{eq:2.2}
\end{equation}
We define the process $W^{(n)}=(W^{(n)}_t)_{0\le t\le T}$ by
\begin{equation*}
W_t^{(n)}:=\frac{1}{\sqrt{n}}\sum_{i=1}^{\lfloor nt\rfloor} 
\xi_i\qquad(0\le
t\le T), 
\end{equation*}
where $\lfloor x\rfloor$ denotes the greatest integer not exceeding $x$. 
The
process $W^{(n)}$ converges weakly to $W$ in the Skorohod space by 
Donsker's
theorem (see, e.g., Billingsley \cite[Theorem 16.1]{B}). We define 
the process $Y^{(n)}=(Y^{(n)}_t)_{0\le t\le T}$ by
\begin{equation*}
Y_t^{(n)}:=\int_{0}^{t}y(\tfrac{\lfloor nt\rfloor}{n},s)dW_s^{(n)}
\qquad(0\le t\le T). 
\end{equation*}
Then it follows that 
\begin{equation*}
Y_t^{(n)}=\frac{1}{\sqrt{n}}\sum_{i=1}^{\lfloor nt\rfloor}  
y(\tfrac{\lfloor
nt\rfloor}{n},\tfrac{i}{n})\xi_i\qquad(0\le t\le T). 
\end{equation*}

Here is the Donsker-type theorem for $Y$.

\begin{thm}
\label{thm:2.1} The process $Y^{(n)}$ converges weakly to $Y$ as 
$n\to\infty$%
.
\end{thm}

\begin{proof}
We first show that the finite-dimensional distributions of 
$Y^{(n)}$ converge to those of $Y$ as $n\to\infty$. 
Thus, for $a_1,\dots,a_d\in\mathbf{R}$ and $t_1,\dots,t_d\in [0,T]$, 
we show that $X^{(n)}$ converges to a normal distribution 
with variance $\mathrm{Var}(X)$, where 
$X^{(n)}:=\sum_{k=1}^{d}a_kY^{(n)}_{t_k}$ and 
$X:=\sum_{k=1}^{d}a_kY_{t_k}$. 
We have
\[
\begin{split}
\mathrm{Var}(X^{(n)})
&=\sum_{k,l=1}^{d}a_ka_l\frac{1}{n}
\sum_{i=1}^{\lfloor n(t_k\wedge t_l)\rfloor}
y(\tfrac{\lfloor nt_k\rfloor}{n},\tfrac{i}{n})
y(\tfrac{\lfloor nt_l\rfloor}{n},\tfrac{i}{n})\\
&=\sum_{k,l=1}^{d}a_ka_l
\int_{0}^{\frac{\lfloor n(t_k\wedge t_l)\rfloor}{n}}
 y(\tfrac{\lfloor nt_k\rfloor}{n},\tfrac{\lfloor ns\rfloor +1}{n})
 y(\tfrac{\lfloor nt_l\rfloor}{n},\tfrac{\lfloor ns\rfloor +1}{n})ds,
\end{split}
\]
where $t\wedge s:=\min(t,s)$. 
The function $(t_1,t_2,u)\mapsto y(t_1,u)y(t_2,u)$ is 
continuous, whence uniformly continuous, on 
the compact set $[0,T]^3$. 
From this and the fact that 
$0\le t-(\lfloor nt\rfloor/n)<1/n$, 
we see that
\begin{equation}
\lim_{n\to\infty}\mathrm{Var}(X^{(n)})
=\sum_{k,l=1}^{d}a_ka_l
\int_{0}^{t_k\wedge t_l}y(t_k,s)
y(t_l,s)ds=\mathrm{Var}(X).
\label{eq:2.3}
\end{equation}

We may assume $\mathrm{Var}(X)>0$. For, otherwise, (\ref{eq:2.3}) implies 
that $X^{(n)}$ converges to $X=0$ in law. 
We put 
$b^{(n)}_i:=\sum_{k=1}^{d}a_ky(\tfrac{\lfloor 
nt_k\rfloor}{n},\tfrac{i}{n})$ 
and 
$X^{(n)}_i:=n^{-1/2}b_i^{(n)}\xi_i$ for $n,i=1,2,\dots$. 
Then we have 
$X^{(n)}=\sum_{i=1}^{\lfloor nT\rfloor}X^{(n)}_i$ for $n=1,2,\dots$. 
We need to show the following Lindeberg's condition: 
for every $\epsilon>0$,
\begin{equation}
\lim_{n\to\infty}
\sum_{i=1}^{\lfloor nT\rfloor}E\left[(X^{(n)}_i)^2
\mathbf{1}_{\{\vert X^{(n)}_i\vert>\epsilon \sigma^{(n)}\}}\right]=0,
\label{eq:2.4}
\end{equation}
where 
$\sigma^{(n)}:=\sqrt{\mathrm{Var}(X^{(n)})}$. 
Choose a positive constant $M$ satisfying 
$\vert b^{(n)}_i\vert \le M$ for $n,i=1,2,\dots$. 
Then since $\vert X^{(n)}_i\vert\le Mn^{-1/2}\vert\xi_i\vert$, we have
\begin{equation*}
\begin{split}
&\sum_{i=1}^{\lfloor nT\rfloor}E\left[(X^{(n)}_i)^2
\mathbf{1}_{\{\vert X^{(n)}_i\vert>\epsilon \sigma^{(n)}\}}\right]
\le 
\sum_{i=1}^{\lfloor nT\rfloor}E\left[(Mn^{-1/2}\xi_i)^2
\mathbf{1}_{\{\vert Mn^{-1/2}\xi_i\vert>\epsilon \sigma^{(n)}\}}\right]\\
&=\sum_{i=1}^{\lfloor nT\rfloor}
M^2n^{-1}
E\left[(\xi_1)^2\mathbf{1}_{\{\vert \xi_1\vert\ge M^{-1}\sigma^{(n)}\sqrt{n}\}}\right]
\le M^2TE\left[(\xi_1)^2
\mathbf{1}_{\{\vert \xi_1\vert\ge M^{-1}\sigma^{(n)}\sqrt{n}\}}\right].
\end{split}
\end{equation*}
We obtain (\ref{eq:2.4}) from this. 
By (\ref{eq:2.4}) and (\ref{eq:2.3}), we can apply 
the central limit theorem (cf.~\cite[Theorem 7.2]{B}), so that 
$X^{(n)}$ converges to $X$ in law, as desired.

Next we show that, for $0\le t_1\le t\le t_2\le T$ and $n=1,2,\dots$,
\begin{equation}
E\left[\vert Y_{t}^{(n)}-Y_{t_1}^{(n)}\vert^2
\vert Y_{t_2}^{(n)}-Y_{t}^{(n)}\vert^2\right]
\le C\vert t_2-t_1\vert^2.
\label{eq:2.5}
\end{equation}
The theorem follows from this and \cite[Theorem 15.6]{B}. 
However, if $t_2-t_1<1/n$, 
then either $t_1$ and $t$ or $t$ and $t_2$ lie in the same subinterval 
$[\frac{m}{n},\frac{m+1}{n})$ 
for some $m$, whence the left hand side of (\ref{eq:2.5}) is zero. 
Therefore we may assume that 
$t_2-t_1\ge 1/n$.

We show that
\begin{equation}
E\left[\vert Y_{t}^{(n)}-Y_{s}^{(n)}\vert^4\right]
\le C\vert t-s\vert^2
\label{eq:2.6}
\end{equation}
for $t$, $s$ and $n$ satisfying
\begin{equation}
0\le s<t\le T,\qquad t-s\ge \tfrac{1}{n}.
\label{eq:2.7}
\end{equation}
This implies (\ref{eq:2.5}) under the condition $t_2-t_1\ge 1/n$ since
\begin{equation*}
\begin{split}
E\left[\vert Y_{t}^{(n)}-Y_{t_1}^{(n)}\vert^2
\vert Y_{t_2}^{(n)}-Y_{t}^{(n)}\vert^2\right]
&\le E\left[\vert Y_{t}^{(n)}-Y_{t_1}^{(n)}\vert^4\right]^{1/2}
     E\left[\vert Y_{t_2}^{(n)}-Y_{t}^{(n)}\vert^4\right]^{1/2}\\
&\le C\vert t-t_1\vert\vert t_2-t\vert\le C\vert t_2-t_1\vert^2.
\end{split}
\end{equation*}
For distinct $i$, $j$, $k$ and $l$, we have
\[
E[(\xi_i)^3\xi_j]=E[(\xi_i)^2\xi_j\xi_k]=E[\xi_i\xi_j\xi_k\xi_l]=0.
\]
Hence, for $t$, $s$ and $n$ satisfying (\ref{eq:2.7}), 
$E[\vert Y_{t}^{(n)}-Y_{s}^{(n)}\vert^4]$ is equal to
\[
\begin{split}
&n^{-2}E\left[\left\{\sum\nolimits_{i=1}^{\lfloor nt\rfloor}
 (y(\tfrac{\lfloor nt\rfloor}{n},\tfrac{i}{n})-
  y(\tfrac{\lfloor ns\rfloor}{n},\tfrac{i}{n}))\xi_i\right\}^4\right]\\
&=E[(\xi_1)^4]n^{-2}\sum_{i=1}^{\lfloor nt\rfloor}
 \{y(\tfrac{\lfloor nt\rfloor}{n},\tfrac{i}{n})
 -y(\tfrac{\lfloor ns\rfloor}{n},\tfrac{i}{n})\}^4\\
&+\frac{6}{n^2}E[(\xi_1)^2]^2\sum_{1\le i<j\le \lfloor nt\rfloor}
 \{y(\tfrac{\lfloor nt\rfloor}{n},\tfrac{i}{n})
  -y(\tfrac{\lfloor ns\rfloor}{n},\tfrac{i}{n})\}^2
 \{y(\tfrac{\lfloor nt\rfloor}{n},\tfrac{j}{n})
  -y(\tfrac{\lfloor ns\rfloor}{n},\tfrac{j}{n})\}^2\\
&=(I_1+I_2)E[(\xi_1)^4]+6(J_1+J_2+J_3)E[(\xi_1)^2]^2,
\end{split}
\]
where
\begin{equation*}
I_1:=n^{-2}\sum_{i=1}^{\lfloor ns\rfloor}
 \{y(\tfrac{\lfloor nt\rfloor}{n},\tfrac{i}{n})
  -y(\tfrac{\lfloor ns\rfloor}{n},\tfrac{i}{n})\}^4,\qquad 
I_2:=n^{-2}\sum_{i=\lfloor ns\rfloor +1}^{\lfloor nt\rfloor}
 y(\tfrac{\lfloor nt\rfloor}{n},\tfrac{i}{n})^4
\end{equation*}
and
\begin{align*}
&J_1:=n^{-2}\sum_{(i,j)\in \Lambda_1}
 \{y(\tfrac{\lfloor nt\rfloor}{n},\tfrac{i}{n})
  -y(\tfrac{\lfloor ns\rfloor}{n},\tfrac{i}{n})\}^2
 \{y(\tfrac{\lfloor nt\rfloor}{n},\tfrac{j}{n})
  -y(\tfrac{\lfloor ns\rfloor}{n},\tfrac{j}{n})\}^2,\\
&J_2:=n^{-2}\sum_{(i,j)\in \Lambda_2}
 \{y(\tfrac{\lfloor nt\rfloor}{n},\tfrac{i}{n})
  -y(\tfrac{\lfloor ns\rfloor}{n},\tfrac{i}{n})\}^2
   y(\tfrac{\lfloor nt\rfloor}{n},\tfrac{j}{n})^2,\\
&J_3:=n^{-2}\sum_{(i,j)\in \Lambda_3}
 y(\tfrac{\lfloor nt\rfloor}{n},\tfrac{i}{n})^2
 y(\tfrac{\lfloor nt\rfloor}{n},\tfrac{j}{n})^2
\end{align*}
with
\begin{align*}
&\Lambda_1:=\{(i,j): 1\le i<j\le \lfloor ns\rfloor\},\\
&\Lambda_2:=\{(i,j): 1\le i\le \lfloor ns\rfloor,\ 
 \lfloor ns\rfloor < j\le \lfloor nt\rfloor\},\\
&\Lambda_2:=\{(i,j): \lfloor ns\rfloor <i<j\le \lfloor nt\rfloor\}.
\end{align*}
By (\ref{eq:2.7}), we have
\[
 \lfloor nt\rfloor -\lfloor ns\rfloor\le 
 nt-ns+1=n(t-s+\tfrac{1}{n})\le 2n(t-s),
\]
so that
\[
\#\Lambda_1\le Cn^2,\quad 
\#\Lambda_2\le Cn^2(t-s),\quad 
\#\Lambda_3\le Cn^2(t-s)^2.
\]
Therefore, using (\ref{eq:2.1}), 
we obtain, for $t$, $s$ and $n$ satisfying (\ref{eq:2.7}),
\begin{align*}
&\vert I_1\vert\le Cn^{-2}\cdot n\cdot(t-s)^4=Cn^{-1}(t-s)^4
\le C(t-s)^5\le C(t-s)^2,\\
&\vert I_2\vert\le Cn^{-2}\cdot n(t-s)=Cn^{-1}(t-s)\le C(t-s)^2,\\
&\vert J_1\vert\le Cn^{-2}\cdot n^2\cdot (t-s)^4=C(t-s)^4\le C(t-s)^2,\\
&\vert J_2\vert\le Cn^{-2}\cdot n^2(t-s)\cdot (t-s)^2=C(t-s)^3\le 
 C(t-s)^2,\\
&\vert J_3\vert\le Cn^{-2}\cdot n^2(t-s)^2=C(t-s)^2.
\end{align*}
Thus (\ref{eq:2.6}) follows.
\end{proof}

Denote by $\Delta X$ and $[X]$ the jump and quadratic variation processes 
of a process $X$, respectively, i.e., 
\begin{equation*}
\Delta X_t:=X_t-\lim_{s\uparrow t}X_s,\qquad [X]_t:=\sum_{s\le
t}^{}\left(\Delta X_s\right)^2. 
\end{equation*}

\begin{thm}
\label{thm:2.2} The process $\Delta Y^{(n)}$ converges to zero in
probability, while $[Y^{(n)}]$ converges to the deterministic process $%
(t)_{0\le t\le T}$ in probability.
\end{thm}

\begin{proof}
From (\ref{eq:2.6}) with (\ref{eq:2.7}), we have
\[
E\left[(\Delta Y_t^{(n)})^4\right]
\le 
E\left[(Y_t^{(n)}-Y_{t-\frac{1}{n}}^{(n)})^4\right]
\le 
Cn^{-2},
\]
so that, as $n\to\infty$,
\[
E\left[\sup_{0\le t\le T}(\Delta Y_t^{(n)})^4\right]
\le E\left[\sum_{t\le T}(\Delta Y_t^{(n)})^4\right]
=\sum_{t\le T}E\left[(\Delta Y_t^{(n)})^4\right]
\le C\frac{nT}{n^2}\to 0.
\]
Thus $\Delta Y^{(n)}$ converges to zero in probability.

We put 
$Z_t^{(n)}:=\int_{0}^{t}z(\tfrac{\lfloor nt\rfloor}{n},s)dW_s^{(n)}$ 
for $0\le t\le T$. 
Then we have $Y_t^{(n)}=W_t^{(n)}-Z_t^{(n)}$, whence
\[
[Y^{(n)}]_t=[W^{(n)}]_t-2\sum_{s\le t}(\Delta W_s^{(n)})
(\Delta Z_s^{(n)})+[Z^{(n)}]_t.
\]

Since $z(u,u)=0$, we have
\[
Z_t^{(n)}-Z_{t-\frac{1}{n}}^{(n)}
=\frac{1}{\sqrt{n}}\sum_{i=1}^{\lfloor nt\rfloor -1}
 \{z(\tfrac{\lfloor nt\rfloor}{n},\tfrac{i}{n})
  -z(\tfrac{\lfloor nt\rfloor -1}{n},\tfrac{i}{n})\}
\xi_i\quad(\mbox{$=0$\quad if $\lfloor nt\rfloor =1$}).
\]
From this and (\ref{eq:2.1}), $E[(\Delta Z_t^{(n)})^2]$ is at most
\[
\begin{split}
E\left[(Z_t^{(n)}-Z_{t-\frac{1}{n}}^{(n)})^2\right]
&=\frac{1}{n}\sum_{i=1}^{\lfloor nt\rfloor -1}
\{z(\tfrac{\lfloor nt\rfloor}{n},\tfrac{i}{n})
 -z(\tfrac{\lfloor nt\rfloor -1}{n},\tfrac{i}{n})\}^2
\le\frac{nT}{n}\cdot \frac{C^2}{n^2}=\frac{C}{n^2}.
\end{split}
\]
Since $[Z^{(n)}]_t$ is increasing, we see that
\begin{equation}
E\left[\sup_{0\le t\le T}[Z^{(n)}]_t\right]
=E\left[[Z^{(n)}]_T\right]=\sum_{t\le T}E\left[(\Delta Z_t^{(n)})^2\right]
\le nT\frac{C}{n^2}=\frac{C}{n}.
\label{eq:2.8}
\end{equation}
Thus $[Z^{(n)}]$ converges to zero in probability.

We have
\begin{equation*}
[W^{(n)}]_t-t=\frac{\lfloor nt\rfloor}{n}-t+
\frac{1}{n}\sum_{i=1}^{\lfloor nt\rfloor}\{(\xi_i)^2-1\}.
\end{equation*}
Let $\epsilon>0$. Then, 
from (\ref{eq:2.2}) and 
Kolmogorov's inequality (see, e.g, Williams \cite[Section 14.6]{W}), we see 
that
\[
\begin{split}
&P\left(\sup_{0\le t\le T}\frac{1}{n}\left\vert
\sum\nolimits_{i=1}^{\lfloor nt\rfloor}
 \{(\xi_i)^2-1\}\right\vert\ge \epsilon\right)
=P\left(\sup_{0\le t\le T}\left\vert\sum\nolimits_{i=1}^{\lfloor nt\rfloor}
 \{(\xi_i)^2-1\}\right\vert\ge n\epsilon\right)\\
&\le 
\frac{1}{\epsilon^2n^2}\sum_{i=1}^{\lfloor nT\rfloor}
E\left[(\xi_i^2-1)^2\right]
\le \frac{nT}{\epsilon^2n^2}E\left[(\xi_1^2-1)^2\right]\to 0\quad(n\to\infty).
\end{split}
\]
From this and the fact that 
$0\le t-(\lfloor nt\rfloor/n)<1/n$, we see that 
$[W^{(n)}]$ converges to the deterministic process $(t)$ in probability.

By Schwarz's inequality, we have
\[
\left\vert\sum\nolimits_{s\le t}(\Delta W_s^{(n)})(\Delta Z_s^{(n)})\right\vert
\le [W^{(n)}]_t^{1/2}[Z^{(n)}]_t^{1/2}
\le [W^{(n)}]_T^{1/2}[Z^{(n)}]_T^{1/2},
\]
whence, by (\ref{eq:2.8}),
\[
\begin{split}
&E\left[\sup_{0\le t\le T}
\left\vert\sum\nolimits_{s\le t}(\Delta W_s^{(n)})(\Delta Z_s^{(n)})
\right\vert\right]
\le 
E\left[[W^{(n)}]_T^{1/2}[Z^{(n)}]_T^{1/2}\right]\\
&\le 
E\left[[W^{(n)}]_T\right]^{1/2}E\left[[Z^{(n)}]_T\right]^{1/2}
\le T^{1/2}\cdot (Cn^{-1})^{1/2}=Cn^{-1/2}.
\end{split}
\]
Thus the process $(\sum_{s\le t}(\Delta W_s^{(n)})(\Delta Z_s^{(n)}))$ 
also converges to zero in probability. 
Combining, we see that $[Y^{(n)}]$ converges to $(t)$ in probability.
\end{proof}


\section{Approximating binary market}
\label{sec:3}

Let $T\in (0,\infty)$ and let $Y$ be as defined in Section \ref{sec:2}. 
We consider the stock price
process $S$ that is governed by the following more general stochastic
differential equaltion than (\ref{eq:1.1}): 
\begin{equation*}
dS_{t}=S_{t}\{b(t)dt+\sigma dY_{t}\}\qquad (0\le t\le T),
\end{equation*}
where $\sigma $ and the initial value $S_{0}$ are positive constants, and 
$%
b(\cdot )$ is a deterministic continuous function on $[0,T]$. The solution 
$S
$ is given by 
\begin{equation*}
S_{t}:=S_{0}\exp \left\{\sigma Y_t+\int_{0}^{t}b(s)ds-\tfrac{1}{2}\sigma
^{2}t\right\}\qquad (0\leq t\leq T).
\end{equation*}
For $n=1,2,\dots$, we consider the process 
$S^{(n)}=(S^{(n)}_t)_{0\le t\le T}$ defined by 
\begin{equation*}
S_{t}^{(n)}:=\prod_{s\leq t}\left\{1+\sigma \Delta Y_{s}^{(n)}+\frac{1}{n}b(%
\tfrac{\lfloor ns\rfloor }{n})\right\}\qquad (0\leq t\leq T),
\end{equation*}
where $Y^{(n)}$ is as in Section 2. The aim of this section is to prove 
that 
$S^{(n)}$ converges weakly to the process $S$.

As in \cite[(10) and (11)]{S}, we put 
\begin{equation*}
Y_t^{(1,n)}:=\sum_{s\le t}\Delta Y^{(n)}_s \mathbf{1}_{\{\vert \Delta
Y_s^{(n)}\vert<\frac{1}{2}\sigma^{-1}\}},\quad Y_t^{(2,n)}:=\sum_{s\le
t}\Delta Y^{(n)}_s \mathbf{1}_{\{\vert \Delta Y_s^{(n)}\vert\ge\frac{1}{2}%
\sigma^{-1}\}}.
\end{equation*}
Then we have 
\begin{align}
&Y^{(n)}_t=Y^{(1,n)}_t+Y^{(2,n)}_t,  \label{eq:3.1} \\
&[Y^{(1,n)}]_t=\sum_{s\le t}(\Delta Y^{(n)}_s)^2 \mathbf{1}_{\{\vert \Delta
Y_s^{(n)}\vert<\frac{1}{2}\sigma^{-1}\}},  \label{eq:3.2} \\
&[Y^{(2,n)}]_t=\sum_{s\le t}(\Delta Y^{(n)}_s)^2 \mathbf{1}_{\{\vert \Delta
Y_s^{(n)}\vert\ge\frac{1}{2}\sigma^{-1}\}},  \label{eq:3.3} \\
&[Y^{(n)}]_t=[Y^{(1,n)}]_t+[Y^{(2,n)}]_t.  \label{eq:3.4}
\end{align}

\begin{lem}
\label{lem:3.1} The process $[Y^{(2,n)}]$ converges to zero in probability,
whence $[Y^{(1,n)}]$ converges to the deterministic process $(t)$ in
probability. The process $Y^{(2,n)}$ converges to zero in probability,
whence $Y^{(1,n)}$ converges weakly to $Y$.
\end{lem}

\begin{proof}
Let $\epsilon>0$. Then, by (\ref{eq:3.3}), we have
\begin{equation*}
P\left(\sup_{0\le t\le T}[Y^{(2,n)}]_t\ge \epsilon\right)
\le P\left(\sup_{0\le t\le T}[Y^{(2,n)}]_t>0\right)
=P\left(\sup_{0\le t\le T}\vert \Delta Y_t^{(n)}\vert\ge 
\tfrac{1}{2}\sigma^{-1}\right).
\end{equation*}
Since the process $\Delta Y^{(n)}$ converges to zero in probability by 
Theorem \ref{thm:2.2}, $[Y^{(2,n)}]$ converges to zero in probability. 
Therefore, by Theorem \ref{thm:2.2} and (\ref{eq:3.4}), $[Y^{(1,n)}]$ 
converges to zero in probability.

In the same way, since
\[
P\left(\sup_{0\le t\le T}\vert Y^{(2,n)}_t\vert\ge \epsilon\right)
\le P\left(\sup_{0\le t\le T}\vert \Delta Y_t^{(n)}\vert\ge \tfrac{1}{2}
\sigma^{-1}\right),
\]
it follows from Theorem \ref{thm:2.2} that 
$Y^{(2,n)}$ converges to zero in probability. Therefore, by 
Theorem \ref{thm:2.1}, (\ref{eq:3.1}) and \cite[Theorem 4.1]{B}, 
$Y^{(1,n)}$ converges weakly to $Y$.
\end{proof}

\begin{thm}
\label{thm:3.2} The process $S^{(n)}$ converges weakly to $S$.
\end{thm}

\begin{proof}
Write 
$S_t^{(n)}=S_t^{(1,n)}S_t^{(2,n)}$, 
where
\begin{align*}
S_t^{(1,n)}&:=\prod_{s\le t}\left\{1+\sigma\Delta Y_s^{(1,n)}+
 \frac{1}{n}b(\tfrac{\lfloor ns\rfloor}{n})\right\} \\
S_t^{(2,n)}&:=\prod_{s\le t}\left\{1+\sigma\Delta Y_s^{(2,n)}\right\}, 
\end{align*}
and the processes $Y^{(i,n)}$ are as above. 
We claim the following: 
(i)\ $S^{(1,n)}$ {\it converges weakly to\/} $S$; 
(ii)\ $S^{(2,n)}$ {\it converges to one in probability\/}.

By \cite[Problem 1, Page 28]{B}, the claim (ii) implies that
$S^{(1,n)}(S^{(2,n)}-1)$ converges to zero in probability. 
Since
\[
S_t^{(n)}=S_t^{(1,n)}(S_t^{(2,n)}-1)+S_t^{(1,n)},
\]
we see from (i) and \cite[Theorem 4.1]{B} that 
$S^{(n)}$ converges weakly to $S$, as desired.

We first prove (ii). Let $\epsilon>0$. Then
\[
P\left(\sup_{0\le t\le T}\vert S_t^{(2,n)}-1\vert\ge \epsilon\right)
\le P\left(\sup_{0\le t\le T}\vert \Delta Y_{t}^{(n)}\vert>
 \tfrac{1}{2}\sigma^{-1}\right).
\]
Since the process $\Delta Y^{(n)}$ converges to zero in probability by 
Theorem \ref{thm:2.2}, $S^{(2,n)}$ converges to one in probability. 
Thus (ii) follows. Next we prove (i). Since the exponential is a continuous 
functional in the Skorohod topology, it is enough to prove that 
$\log S^{(1,n)}$ converges weakly to the process 
$(\sigma Y_t+\int_{0}^{t}b(s)ds-\frac{1}{2}\sigma^2t)$. 
Notice that 
$\vert \sigma\Delta Y_t^{(1,n)}\vert 
 +\frac{1}{n}\vert b(\tfrac{\lfloor nt\rfloor}{n})\vert<\frac{3}{4}$ 
for sufficiently large $n$ and $t\in [0,T]$, whence the logarithm 
$\log S^{(1,n)}$ is well defined for such $n$.

We have
\[
\log (1+x)=x-\frac{1}{2}x^2+r(x)x^3\qquad
(\vert x\vert<1),
\]
where $r(x)$ is a bounded function on $\vert x\vert\le \frac{3}{4}$. 
Hence
\begin{align*}
\log S_t^{(1,n)}
&=\sum_{s\le t}\left\{\sigma\Delta Y_s^{(1,n)}
+\frac{1}{n}b(\tfrac{\lfloor ns\rfloor}{n})
-\frac{1}{2}\left(\sigma\Delta Y_s^{(1,n)}+
\frac{1}{n}b(\tfrac{\lfloor ns\rfloor}{n})\right)^2 \right.\\
&\qquad\qquad\quad\left. +r\left(\sigma\Delta Y_s^{(1,n)}
  +\frac{1}{n}b(\tfrac{\lfloor ns\rfloor}{n})\right)\cdot
 \left(\sigma\Delta Y_s^{(1,n)}+
  \frac{1}{n}b(\tfrac{\lfloor ns\rfloor}{n})\right)^3\right\}\\
&=\sigma Y_t^{(1,n)}+ \sum_{s\le t}\frac{1}{n}
 b(\tfrac{\lfloor ns\rfloor}{n}) -\frac{1}{2}\Phi^{(n)}_t +\Psi^{(n)}_t, 
\end{align*}
where
\begin{align*}
\Phi^{(n)}_t&:= \sum_{s\le t}\left(\frac{1}{n}b(\tfrac{\lfloor ns\rfloor}{n})+
 \sigma\Delta Y_s^{(1,n)}\right)^2, \\
\Psi^{(n)}_t&:= \sum_{s\le t}r\left(\sigma\Delta Y_s^{(1,n)}
  +\frac{1}{n}b(\tfrac{\lfloor ns\rfloor}{n})\right)\cdot
 \left(\sigma\Delta Y_s^{(1,n)}+
  \frac{1}{n}b(\tfrac{\lfloor ns\rfloor}{n})\right)^3.  
\end{align*}

We have $\Phi^{(n)}_t=n^{-2}\sum_{s\le t}b(\tfrac{\lfloor ns\rfloor}{n})^2
 +2\sigma\Gamma^{(N)}_t+\sigma^2[Y^{(1,n)}]_t$, where
\begin{equation*}
\Gamma^{(n)}_t:=\sum\nolimits_{s\le t}\frac{1}{n}
 b(\tfrac{\lfloor ns\rfloor}{n})\Delta Y_s^{(1,n)}.
\end{equation*}
Since $b(\cdot)$ is bounded, the first term 
$n^{-2}\sum_{s\le t}b(\tfrac{\lfloor ns\rfloor}{n})^2$ goes to 0 
as $n\to\infty$. By Lemma \ref{lem:3.1}, 
the third term $\sigma^2[Y^{(1,n)}]$ converges to $(\sigma^2t)$ in 
probability. 
As for the second term, it holds that
\begin{equation*}
\sup_{0\le t\le T}\left\vert\Gamma^{(n)}_t\right\vert
\le C\sup_{s\le T}\vert \Delta Y_s^{(1,n)}\vert
\le C\le \vert\Delta Y_t^{(n)}\vert.
\end{equation*}
Since $\Delta Y^{(n)}$ converges to zero in probability by Theorem 
\ref{thm:2.2}, so does $\Gamma^{(n)}$. Thus the process $(\Phi_t)$ 
converges to $(\sigma^2t)$. 
Since
\begin{equation*}
\sup_{0\le t\le T}\Psi_t 
\le C\left(\frac{1}{n}+ \sup_{s\le T}\vert \Delta Y_s^{(1,n)}\vert
 \right)\Phi_T,
\end{equation*}
we see that the process $(\Psi_t)$ converges to zero in probability. 
Using these fact as well as Lemma \ref{lem:3.1} and 
\cite[Theorem 4.1]{B}, we see that $\log S^{(1,n)}$ converges weakly 
to $(\sigma Y_t +\int_0^t b(s)ds -\frac{1}{2}\sigma^2t)$.
\end{proof}

If we take the i.i.d.~random variables $\{\xi_i\}$ so that 
\begin{equation}
P(\xi_1=1)=P(\xi_1=-1)=1/2,  \label{eq:3.5}
\end{equation}
then we obtain the desired approximating binary market model.


\section{Arbitrage opportunities in the binary market}

\label{sec:4}

In this section, we study the arbitrage opportunities in the approximating
binary market model with memory constructed in Section \ref{sec:3}. For
simplicity, we assume that the function $b(\cdot)$ is a real constant as in 
(%
\ref{eq:1.1}).

Let $N\in\mathbf{N}$, $r, b\in\mathbf{R}$, and $\sigma\in (0,\infty)$. The
number $N$ corresponds to $n$ in Sections \ref{sec:2} and \ref{sec:3}. Let
the function $y(t,u)$ be as in Section \ref{sec:2}. We define 
\begin{equation*}
r^{(N)}:=\frac{r}{N},\qquad b^{(N)}:=\frac{b}{N}. 
\end{equation*}
The $\lfloor NT\rfloor$-period market $\mathcal{M}^{(N)}$ 
consists
of a share of the money market with price process $(B^{(N)}_n)_{n=0,1,\dots,%
\lfloor NT\rfloor}$ and a stock with price process $(S^{(N)}_n)_{n=0,1,\dots,%
\lfloor NT\rfloor}$. The prices are governed respectively by 
\begin{gather*}
B^{(N)}_n=B^{(N)}_{n-1}(1+r^{(N)})\quad (n=1,\dots,\lfloor NT\rfloor),
\quad B^{(N)}_0=1, \\
S^{(N)}_n=S^{(N)}_{n-1}(1+ b^{(N)}+ X^{(N)}_n)\quad (n=1,\dots,\lfloor NT\rfloor),
\quad S^{(N)}_0=s_0,
\end{gather*}
where $s_0$ is a positive constant, 
\begin{equation*}
X^{(N)}_n:=\sigma\Delta Y_{\frac{n}{N}}^{(N)}  =\frac{\sigma}{\sqrt{N}}%
\sum_{i=1}^n \left\{y(\tfrac{n}{N},\tfrac{i}{N})-  
y(\tfrac{n-1}{N},\tfrac{i%
}{N})\right\}\xi_i 
\end{equation*}
and $\{\xi_i\}$ are i.i.d.~random variables satisfying (\ref{eq:3.5}). 
Theorem \ref{thm:3.2} implies that the binary market model $\mathcal{M}%
^{(N)}$ approximates the continuous-time market model with bond price
process $(e^{rt})$ and stock price process $S$ in (\ref{eq:1.1}).

Given the values of $\xi_1,\dots,\xi_{n-1}$, the random variable $X^{(N)}_n$ 
takes
the following two possible values $u_n$ and $d_n$: $d_1=-\sigma/\sqrt{N},\;
u_1=\sigma/\sqrt{N}$, and for $n=2,\dots, N$, 
\begin{align*}
d_n&\equiv d_n(\xi_1,\dots,\xi_{n-1})= \frac{\sigma}{\sqrt{N}}%
\sum_{i=1}^{n-1} \left\{y(\tfrac{n}{N},\tfrac{i}{N})- y(\tfrac{n-1}{N},%
\tfrac{i}{N})\right\}\xi_i -\frac{\sigma}{\sqrt{N}}, \\
u_n&\equiv u_n(\xi_1,\dots,\xi_{n-1})= \frac{\sigma}{\sqrt{N}}%
\sum_{i=1}^{n-1} \left\{y(\tfrac{n}{N},\tfrac{i}{N})- y(\tfrac{n-1}{N},%
\tfrac{i}{N})\right\}\xi_i +\frac{\sigma}{\sqrt{N}}.
\end{align*}

We investigate the arbitrage opportunities in $\mathcal{M}^{(N)}$. Let $C$
be a positive constant satisfying 
\begin{equation}  \label{eq:4.1}
|y(t,u)- y(s,u)|\le C |t-s| \qquad (0\le t,s,u\le T).
\end{equation}

\begin{thm}
\label{thm:4.1} Suppose that $T<1/C$. Then there exists an integer $N_0$
such that for each $N\ge N_0$, the market $\mathcal{M}^{(N)}$ is
arbitrage-free.
\end{thm}

\begin{proof}
From the condition $TC<1$, we have an integer $N_0$ satisfying 
\begin{equation}
\label{eq:4.2}
 \frac{b}{N}-\frac{\sigma}{\sqrt{N}}(TC +1)
  > -1,  \quad
 |r-b|< \sqrt{N}(1-TC)\sigma\quad (N\ge N_0).
\end{equation}
By (\ref{eq:4.1}), we have, for $n=1,\dots,\lfloor NT\rfloor$,
\begin{align*}
 \min_{\xi\in\left\{-1,1\right\}^{n-1}}d_n(\xi)
 &= -\frac{\sigma}{\sqrt{N}}\sum_{i=1}^{n-1}
   \vert y(\tfrac{n}{N},\tfrac{i}{N})- 
    y(\tfrac{n-1}{N},\tfrac{i}{N})\vert -\frac{\sigma}{\sqrt{N}} \\
 &\ge -\frac{\sigma}{\sqrt{N}}\left(\frac{n-1}{N}C+1\right)
 \ge -\frac{\sigma}{\sqrt{N}}\left(TC+1\right).   
\end{align*}
This and (\ref{eq:4.2}) yield, 
for $N\ge N_0$ and $n=1,\dots,\lfloor NT\rfloor$,
\begin{equation*}
 b^{(N)}+ X^{(N)}_n\ge 
  \frac{b}{N} + \min_{\xi\in\left\{-1,1\right\}^{n-1}}d_n(\xi) 
  > -1,
\end{equation*}
whence $S_n>0$.

We show that $\mathcal{M}^{(N)}$ is arbitrage-free for $N\ge N_0$. 
By Dzhaparidze \cite[Proposition 6.1.2]{D}, $\mathcal{M}^{(N)}$ 
is free from arbitrage opportunities if and only if 
\begin{equation}
 d_n < r^{(N)} - b^{(N)} < u_n \qquad (n=1,\dots,\lfloor NT\rfloor).
\label{eq:4.3} 
\end{equation}
However, we have
\begin{align*}
 \max_{\xi\in\left\{-1,1\right\}^{n-1}} d_n(\xi)
 &= \frac{\sigma}{\sqrt{N}}\sum_{i=1}^{n-1}
    \vert y(\tfrac{n}{N},\tfrac{i}{N})- 
    y(\tfrac{n-1}{N},\tfrac{i}{N})\vert -\frac{\sigma}{\sqrt{N}} \\
 &\le -\frac{\sigma}{\sqrt{N}}
       \left(1-\frac{n-1}{N}C\right)\le -\frac{\sigma}{\sqrt{N}}
       \left(1-TC\right),
\end{align*} 
and
\begin{align*}
 \min_{\xi\in\left\{-1,1\right\}^{n-1}} u_n(\xi)
 &= -\frac{\sigma}{\sqrt{N}}\sum_{i=1}^{n-1}
    \vert y(\tfrac{n}{N},\tfrac{i}{N})- 
    y(\tfrac{n-1}{N},\tfrac{i}{N})\vert +\frac{\sigma}{\sqrt{N}} \\
 &\ge \frac{\sigma}{\sqrt{N}}
       \left(1-\frac{n-1}{N}C\right)\ge \frac{\sigma}{\sqrt{N}}
       \left(1-TC\right).  
\end{align*}
Thus, by (\ref{eq:4.2}), (\ref{eq:4.3}) holds for $N\ge N_0$.
\end{proof}

By Theorem \ref{thm:4.1}, the market $\mathcal{M}^{(N)}$ is arbitrage-free
for $T$ small enough and $N$ large enough. However, in general, the market 
$%
\mathcal{M}^{(N)}$ may admit arbitrage opportunities, as we see below.

Suppose that there exists a positive constant $C$ such that $l(s,u)\geq C$
for $0\leq u<s\leq T$. Let $T>1/C$. We assume that $r\leq b$. Then, $%
d_{\lfloor NT\rfloor }(-1,\dots ,-1)$ is 
\begin{equation*}
\frac{\sigma }{\sqrt{N}}\sum_{i=1}^{\lfloor NT\rfloor -1}\int_{\tfrac{%
\lfloor NT\rfloor -1}{N}}^{\tfrac{\lfloor NT\rfloor}{N}}l(s,\tfrac{i}{N})ds-\frac{\sigma 
}{%
\sqrt{N}}>\frac{\sigma }{\sqrt{N}}\left( \frac{C(\lfloor NT\rfloor -1)}{N}%
-1\right) .
\end{equation*}
Since $TC>1$, it follows that $d_{\lfloor NT\rfloor }(-1,\dots
,-1)>r_{N}-b_{N}$ or 
\begin{equation*}
S_{\lfloor NT\rfloor }>(1+r_{N})S_{\lfloor NT\rfloor -1}
\end{equation*}
for $N$ large enough. Therefore, if the value of $(\xi _{1},\dots ,\xi
_{\lfloor NT\rfloor -1})$ turns out to be $(-1,\dots ,-1)$, then we have an
arbitrage opportunity: we may buy stocks at time $\lfloor NT\rfloor -1$
using money obtained by shortselling bonds. In a similar fashion, we can
show that if $T>1/C$, $r<b$ and $N$ is large enough, then 
the value $(1,\dots ,1)$ 
of 
$(\xi _{1},\dots ,\xi _{\lfloor NT\rfloor -1})$ gives an arbitrage
opportunity.

Put 
\begin{equation*}
P_N= P\left(\bigcup\nolimits_{n=1}^{\lfloor NT\rfloor} \left\{d_n< r^{(N)}-b^{(N)}<
u_n\right\}^{c}\right). 
\end{equation*}
As we see in the proof of Theorem \ref{thm:4.1}, the binary market 
$\mathcal{%
M}^{(N)}$ is arbitrage-free if and only if $P_N=0$. The next theorem gives
the rate at which the arbitrage probability $P_N$ tends to zero as $%
N\to\infty$.

\begin{thm}
\label{thm:4.2}  There exists a positive constant $C^{\prime}=C^{\prime}_T$
such that, for each $\alpha\in (0,1)$, we have $N(\alpha)\in\mathbf{N}$
satisfying 
\begin{equation*}
P_N\le \frac{C^{\prime}}{N^{\alpha}} \qquad (N\ge N(\alpha)). 
\end{equation*}
\end{thm}

\begin{proof}
Set $\beta:=(\alpha+1)/2$, and choose $N(\alpha)\in\mathbf{N}$ 
so large that
\begin{equation}
N^{\beta/2}C\sqrt{T}< 
  \sqrt{N}-|(r-b)/\sigma|,\qquad 
  N^{\beta/2}>4\qquad(N\ge N(\alpha)).
\label{eq:4.4}
\end{equation}
Then we have $d_1<r^{(N)}-b^{(N)}<u_1$. 
For $N\ge N(\alpha)$ and 
$n=2,\dots,\lfloor NT\rfloor$, we put $\lambda :=N^{\beta/2}$ 
and
\begin{equation*}
s_{n-1}:=\left[N\sum_{i=1}^{n-1}\left\{y(\tfrac{n}{N},\tfrac{i}{N})
  -y(\tfrac{n-1}{N},\tfrac{i}{N})\right\}^2\right]^{1/2}, \qquad
M_{n-1}:=\max_{1\le m\le n-1}
\left\vert\sum\nolimits_{i=1}^{m}\eta_i\right\vert,
\end{equation*}
where $\eta_i:=\sqrt{N}\left\{y(\tfrac{n}{N},\tfrac{i}{N})
  -y(\tfrac{n-1}{N},\tfrac{i}{N})\right\}\xi_i$ for $i=1,2,\dots$. 
By (\ref{eq:4.1}), we have $s_{n-1}\le C\sqrt{T}$. This and (\ref{eq:4.4}) 
imply that
\begin{align*}
P\left(\frac{r-b}{N}\le d_n\right)&\le 
P\left(\frac{r-b}{\sigma}+ \sqrt{N}\le M_{n-1}\right)\le 
P(M_{n-1}\ge \lambda C\sqrt{T}) \\
&\le P(M_{n-1}\ge \lambda s_{n-1}). 
\end{align*}
Similarly we have
\[
 P\left(u_n\le \frac{r-b}{N}\right)\le P(M_{n-1}\ge \lambda s_{n-1}).
\]
Since $\frac{1}{4}\lambda >1$ and 
\begin{equation*}
\max_{1\le i\le n-1}| \eta_i|
=\max_{1\le i\le n-1}| \sqrt{N}\{y(\tfrac{n}{N},\tfrac{i}{N})
  -y(\tfrac{n-1}{N},\tfrac{i}{N})\}|\le s_{n-1},
\end{equation*}
it follows from \cite[(12.16), Page 89]{B} that 
\[
 P(M_{n-1}\ge \lambda s_{n-1})\le \frac{C_0}{\lambda^4}
\]
for some constant $C_0>0$ independent of $N$ and $n$ 
(notice that $\eta_i$ here corresponds to $\xi_i$ in 
\cite[(12.16), Page 89]{B}). 
Hence, $P_N$ is at most
\begin{align*}
\sum_{n=2}^{\lfloor NT\rfloor}
 \left\{P\left(\frac{r-b}{N}\le d_n\right)+
 P\left(u_n\le\frac{r-b}{N}\right)\right\}
 \le \frac{2\lfloor NT\rfloor C_0}
  {N^{2\beta}}\le\frac{2TC_0}{N^{\alpha}}.
\end{align*}
Thus the theorem follows.
\end{proof}

\end{document}